\documentclass{amsart}
\usepackage{amscd,amssymb,hyperref}
\usepackage[all]{xy}
\theoremstyle{plain}
\newtheorem{prop}[subsection]{Proposition}
\newtheorem{thm}[subsection]{Theorem}
\newtheorem{lem}[subsection]{Lemma}
\newtheorem{cor}[subsection]{Corollary}

\theoremstyle{definition}
\newtheorem{defn}[subsection]{Definition}

\numberwithin{equation}{section}

\begin{document}

\title{Configurations, and parallelograms associated to centers of mass}

\author[F.~R.~Cohen]{F. R.~Cohen$^{*}$}
\address{Department of Mathematics, University of Rochester,
Rochester, NY 14627}
\email{\href{mailto:cohf@math.rochester.edu}{cohf@math.rochester.edu}}
\thanks{{$^{*}$}Partially supported by the NSF Grant No. DMS-0072173 and CNRS-NSF Grant No. 17149}

\author[Y.~Kamiyama]{Y.~Kamiyama$^{\*}$}
\address {Department of Mathematics, University of the Ryukyus,
Nishihara-Cho, Okinawa 903-0213, Japan }
\email{\href{mailto:kamiyama@sci.u-ryukyu.ac.jp}{kamiyama@sci.u-ryukyu.ac.jp}}


\begin{abstract}
The purpose of this article is to

\begin{enumerate}
    \item define $M(t,k)$ the $t$-fold center of mass arrangement
    for $k$ points in the plane,
    \item  give elementary properties of $M(t,k)$ and
    \item give consequences concerning the space $M(2,k)$
of $k$ distinct points in the plane, no four of which are the
vertices of a parallelogram.
\end{enumerate}

The main result proven in this article is that the classical
unordered configuration of $k$ points in the plane is not a retract
up to homotopy of the space of $k$ unordered distinct points in the
plane, no four of which are the vertices of a parallelogram. The
proof below is homotopy theoretic without an explicit computation of
the homology of these spaces.

In addition, a second, speculative part of this article arises from
the failure of these methods in the case of odd primes $p$. This
failure gives rise to a candidate for the localization at odd primes
$p$ of the double loop space of an odd sphere obtained from the
$p$-fold center of mass arrangement. Potential consequences are
listed.
\end{abstract}

\maketitle

\section{Introduction and statement of results}

Fix integers $k$ and $t$. The $t$-fold center of mass arrangement
$M(t,k)$ for integers $t$ with $ k \geq t \geq 1 $ is defined as the
subspace of the $k$-fold product $\mathbb C^k$ given by ordered
k-tuples of points $(x_1,\cdots,x_k)$ such that the centroids of any
set of $t$ elements in the underlying set $\{x_1,\cdots,x_k\}$
$$\sigma_t(x_{i_1},x_{i_2},\cdots,x_{i_t}) = (1/t)(
x_{i_1}+x_{i_2}+\cdots+x_{i_t})$$ are distinct for all distinct
subsets $\{x_{i_1},x_{i_2},\cdots,x_{i_t}\}$, and
$\{x_{j_1},x_{j_2},\cdots,x_{j_t}\}$. In particular, $M(t,k)$ is the
complement of the union of the hyperplanes specified by
$$\sigma_t(x_{i_1},x_{i_2},\cdots,x_{i_t})-\sigma_t(x_{j_1},x_{j_2},\cdots,x_{j_t})
= 0$$ for all pairs of unequal sets $S_I=
\{x_{i_1},x_{i_2},\cdots,x_{i_t}\}$, and $S_J=
\{x_{j_1},x_{j_2},\cdots ,x_{j_t}\}$. Write $$|S_J|$$ for the
cardinality of the set $S_J$. In case $ k < t$, define $M(t,k)$ to
be the Fadell-Neuwirth configuration space $Conf(\mathbb C,k)$ of
ordered $k$ tuples of distinct points in $\mathbb C$ \cite{FN}.

Finite unions of complex hyperplanes in complex $k$-space are known
as complex hyperplane arrangements in \cite{OrlikTerao}. The space
$M(t,k)$ is a complement of a complex hyperplane arrangement.
Consider an equivalent formulation of $M(t,k)$ as the complement of
the variety $V(t,k)$ of ordered $k$-tuples $(x_1,...,x_k)$ defined
by the equation $$\prod_{S_I \neq S_J,|S_I|=|S_J|=t}
([x_{i_1}+x_{i_2}+\cdots +x_{i_t}]-[x_{j_1}+ x_{j_2}+\cdots
+x_{j_t}]) = 0$$ with $$M(t,k) = \mathbb {C}^k - V(t,k).$$

Modifications of the $M(t,k)$, $M'(t,k)$, are defined as follows:
$$M'(t,k) = \cap_{1 \leq s \leq t}M(s,k).$$ Thus $M'(t,k)$ is the
complement of the variety $W(t,k)$ of ordered $k$-tuples
$(x_1,...,x_k)$ defined by the equation $$\prod_{S_I \neq S_J,1<q =
|S_I|= |S_J|\leq t} ([x_{i_1}+x_{i_2}+ \cdots +x_{i_q}]-[x_{j_1}+
x_{j_2}+ \cdots +x_{j_q}]) = 0$$ with $$M'(t,k) = \mathbb C^k -
W(t,k).$$ Similarly, if $ k < t$, define $M'(t,k)$ to be
$Conf(\mathbb C,k)$.

In addition, there are natural inclusions
\[
\begin{CD}
M'(t,k) @>{}>> M(t,k) @>{}>> Conf(\mathbb C,k).\\
\end{CD}
\]  These inclusions are equivariant with
respect to the natural action of the symmetric group on $k$ letters,
$\Sigma_k$.

Consider the $t$-fold symmetric product $\mathbb C^t/ \Sigma_t$, and
notice that there is a map $$\chi_t: Conf(\mathbb C,k) \to (\mathbb
C^t/\Sigma_t)^{k \choose t}$$ gotten by choosing all $t$-element
subsets out of a set of cardinality $k$ with a fixed ordering of the
subsets. The map $\chi_t$ is given on the level of points by the
formula $$\chi_t(z_1,z_2,\ldots, z_k) = \prod_{i_1< i_2 < \ldots <
i_t} [ z_{i_1},z_{i_2},\ldots, z_{i_t}]$$ for which the points $[
z_{i_1},z_{i_2},\ldots, z_{i_t}]$ in $\mathbb C^t/\Sigma_t$ are
ordered in the product left lexicographically by indices and over
all subsets of cardinality $t$ in the set $\{z_1,z_2,\ldots, z_k\}$.

Notice that the map $\chi_t: Conf(\mathbb C,k) \to (\mathbb
C^t/\Sigma_t)^{k \choose t}$ takes values in the configuration space
$ Conf(\mathbb C^t/\Sigma_t, {k \choose t})$. Thus in what follows
below $\chi_t$ will be regarded as a map $$\chi_t: Conf(\mathbb C,k)
\to Conf(\mathbb C^t/\Sigma_t, {k \choose t}).$$

Addition of complex numbers provides a map
$$\oplus_t : \mathbb C^t/\Sigma_t \to \mathbb C $$ with
$$\oplus_t ([z_1,\ldots,z_t])= z_1+ \cdots +z_t.$$
There is an induced map $$\Theta_t: Conf(\mathbb C,k) \to\ {\mathbb
C}^{k \choose t}$$ given by the composite
\[
\begin{CD}
 Conf( \mathbb C,k) @>{\chi_t}>> (\mathbb C^t/\Sigma_t)^ {k \choose t} @>{
 (\oplus_t)^{k \choose t} }>> {\mathbb C}^{k \choose t}.     \\
\end{CD}
\] Thus $$\Theta_t(z_1,z_2,\ldots, z_k) = \prod_{i_1< i_2 < \cdots
< i_t} ( z_{i_1}+ z_{i_2}+ \cdots + z_{i_t})$$ in ${\mathbb C}^{k
\choose t}$.

The next proposition, a useful observation, is recorded next where
$$j:Conf(\mathbb C, {k \choose t}) \to {\mathbb C}^{k \choose t}$$ is
the natural inclusion. This observation is the starting point of the
results here, and provides the basic motivation for considering the
center of mass arrangement .
\begin{prop}\label{prop:proposition cartesion}
The following diagram is a pull-back (a cartesion diagram):
\[
\begin{CD}$$
M(t,k)       @>>>               Conf(\mathbb C, {k \choose t})        \\
@VVV                    @VV\text{j}V                    \\
Conf(\mathbb C,k) @>{\Theta_t}>> {\mathbb C}^{k \choose t}\\
$$\end{CD}
\]
\end{prop}

Notice that $M(2,k)$ is the space of ordered $k$-tuples of distinct
points such that no four of the points are the vertices of a
possibly degenerate parallelogram. Consider the natural inclusion
$M(2,k) \to\ Conf(\mathbb C,k)$ modulo the action of $\Sigma_k$ the
symmetric group on $k$ letters
$$i(2,k):M(2,k)/\Sigma_k \to\ Conf(\mathbb C,k)/\Sigma_k.$$

One question is whether there is a cross-section up to homotopy, or
even a $2$-local stable cross-section up to homotopy for this
inclusion. This last question concerns plane geometry and whether
the configuration space of distinct unordered $k$-tuples of points
in the plane can be deformed to the subspace of points, no four of
which are the vertices of a parallelogram.

\begin{thm} \label{thm:no retraction}
If $k \geq 4$, the natural map
$$i(2,k):M(2,k)/\Sigma_k \to Conf(\mathbb C,k)/\Sigma_k $$ does not admit
a surjection in mod-$2$ homology, and thus does not admit a
cross-section (or a stable $2$-local cross-section) up to homotopy.
\end{thm}

The proof, homotopy theoretic without a specific computation of the
homology of these spaces, gives features of the topology of double
loop spaces which forces the maps $i(2,k):M(2,k)/\Sigma_k \to
Conf(\mathbb C,k)/\Sigma_k $ for $k \geq 4$ to fail to be
epimorphisms in mod-$2$ homology.  The analogous methods applied to
the natural inclusion $$i(p,k): M(p,k)/\Sigma_k \to Conf(\mathbb
C,k)/\Sigma_k $$ for $p$ an odd prime fail to produce a non-trivial
obstruction to the existence of a stable $p$-local section. Hence a
problem unsolved here is whether $i(p,k)$ admits a stable $p$-local
cross-section. The failure of the methods here in case $p$ is an odd
prime leads to the speculation in section $2$ here concerning the
localization of the double loop space of a sphere at an odd prime
$p$.

Further properties of these arrangements are noted next. The natural
``stabilization" map for configuration spaces fails to preserve the
spaces $M(t,k)$. However, there are stabilization maps for the
modified center of mass arrangements $$S: M'(t,k) \to\ M'(t,k+1)$$
defined by $$S(x_1,\ldots,x_k) = (x_1,\ldots,x_k,\vec z)$$ where
$\vec z$ is the vector $(L,0)$ with $ L = 2t( 1 + max_{k \geq i \geq
1} ||x_i||)$. Notice that $S$ takes values in $M'(t,k+1)$, but that
the analogous map out of $ M(t,k)$ takes values in $Conf(\mathbb
C,k)$, but not in the subspace $M(t,k+1)$.

The next result follows directly from \cite{C2,CMT}.
\begin{thm} \label{thm: splitting}
The map $$S: M'(t,k) \to\ M'(t,k+1)$$ extends to a map
$$S_*: M'(t,k)\times_{\Sigma_k}Y^k \to\
M'(t,k+1)\times_{\Sigma_{k+1}}Y^{k+1}$$ which admits a stable left
inverse for any path-connected CW-complex $Y$, and thus induces a
split monomorphism in homology with any field coefficients.
\end{thm}

\begin{cor}\label{cor:sign.rep.splitting} The map $S: M'(t,k)/\Sigma_k  \to\
M'(t,k+1)/ \Sigma_{k+1}$ induces a split monomorphism in homology
with coefficients in any graded permutation representation of
$\Sigma_k$, and thus by specialization to either coefficients given
by the trivial representation or the sign representation.
\end{cor}

Connections to homotopy theory as well as the motivation for
considering the spaces $M(t,k)$ and the map
$$\chi_t: Conf(\mathbb C,k) \to (\mathbb C^t/\Sigma_t)^{k \choose
t}$$ defined earlier in this section are given next. These
connections arise from stable homotopy equivalences
$$H:\Omega^2 \Sigma^2(X) \to \vee_{0 \leq
k}D_k(\Omega^2\Sigma^2(X))$$ in case $X$ is a path-connected
$CW$-complex originally proven in \cite{Snaith} and subsequently in
\cite{CMT,C2} for which $D_k(\Omega^2\Sigma^2(X))$ is defined in
section \ref{speculation} here.

This stable homotopy equivalence is obtained by adding maps given by
$$h_k:\Omega^2 \Sigma^2(X) \to \Omega^{2k}
\Sigma^{2k}D_k(\Omega^2\Sigma^2(X))$$ as observed in the appendix of
\cite{C2}. These maps do not compress through $$\Omega^{2k-1}
\Sigma^{2k-1}D_k(\Omega^2\Sigma^2(X))$$ in case $k = 2^t$, and
spaces are localized at the prime $2$ \cite{CM}.

Specialize $h_k$ to $k=p$ an odd prime and $X = S^{2n-1}$. The
spaces $M(p,k)$ and $M'(p,k)$ as well as the map $\chi_t:
Conf(\mathbb C,k) \to (\mathbb C^t/\Sigma_t)^{k \choose t}$ are
introduced here in order to attempt to compress the maps
$$h_p:\Omega^2 \Sigma^2(S^{2n-1}) \to \Omega^{2p}
\Sigma^{2p}D_p(\Omega^2\Sigma^2(S^{2n-1}))$$ through some choice of
map $$\bar h_p:\Omega^2 \Sigma^2(S^{2n-1}) \to \Omega^{2}
\Sigma^{2}D_p(\Omega^2\Sigma^2(S^{2n-1})).$$

The map $h_p$ as given in \cite{C2,CMT} is induced on the level of
certain combinatorial models by the composite
\[
\begin{CD}
 Conf( \mathbb C,k) @>{\chi_p}>> Conf( \mathbb C^p/\Sigma_p ,{k \choose p} ) @>{inclusion}>>
 (\mathbb C^p/\Sigma_p)^ {k \choose p}.     \\
\end{CD}
\] A space $M_p(\mathbb C, X)$ together with a map
$$I_p:M_p(\mathbb C, X) \to \Omega^2 \Sigma^2(X)$$ will be defined in
section \ref{speculation} in which configuration spaces
$Conf(\mathbb C,k)$ used in combinatorial models of $\Omega^2
\Sigma^2(X)$ are replaced by the spaces $M(p,k)$. Furthermore, there
are continuous maps $$ h_p: M_p(\mathbb C, X) \to
\Omega^{2}\Sigma^{2}(D_p(\Omega^2\Sigma^2(X))).$$

It is natural to compare the homotopy types of $\Omega^2S^{2n+1}$
and $M_p(\mathbb C, S^{2n-1})$ after localization at an odd prime
$p$ by the following theorem in which $$E:\Sigma^2(Y) \to
\Omega^{2p-2}\Sigma^{2p}(Y)$$ denotes the classical suspension map.
\begin{thm} \label{thm:comparison}
There is a commutative diagram
\[
\begin{CD}
M_p(\mathbb C, X)  @>{h_p}>>  \Omega^{2}\Sigma^{2}(D_p(\Omega^2\Sigma^2(X))) \\
 @VV{I_p}V          @VV{\Omega^{2}(E)}V \\
C(\mathbb C, X) @>{h_p}>>
\Omega^{2p}\Sigma^{2p}(D_p(\Omega^2\Sigma^2(X))).
\end{CD}
\]  Thus if the map
$$I_p:M_p(\mathbb C, S^{2n-1})  \to \Omega^2S^{2n+1}$$ is a
$p$-local equivalence, then there is a $p$-local map
$$ \bar h_p: \Omega^2S^{2n+1} \to
\Omega^{2}\Sigma^{2}(D_p(\Omega^2\Sigma^2(S^{2n-1})))$$ which is a
compression of the map $h_p:\Omega^2 \Sigma^2(S^{2n-1}) \to
\Omega^{2p} \Sigma^{2p}D_p(\Omega^2\Sigma^2(S^{2n-1}))$ and which
induces an isomorphism on $H_{2np-2}(-;\mathbb F_p)$.
\end{thm}

Some consequences of the existence of $\bar h_p$ are discussed in
section \ref{speculation} here. These consequences suggest that it
would be interesting to understand the behavior of the natural map
$$I_p: M_p(\mathbb C, S^{2n-1}) \to C(\mathbb C, S^{2n-1})$$ on the
level of mod-$p$ homology.\vskip .2in

 {\bf TABLE OF CONTENTS}
\begin{enumerate}
    \item Introduction
    \item Speculation concerning the localization of the
    double loop space of a sphere at an odd prime $p$, and applications
    \item Sketch of Proposition \ref{prop:proposition cartesion}
    \item Calculations at the prime 2, and the proof of Theorem \ref{thm:no retraction}
    \item Sketch of Theorem \ref{thm: splitting} and Corollary \ref{cor:sign.rep.splitting}
    \item Sketch of Theorem  \ref{thm:comparison}
\end{enumerate}

The authors would like to congratulate Huynh Mui on this happy
occasion of his $60$-th birthday. The work here is inspired by Mui's
mathematical work on extended power constructions as well as his
interest in configuration spaces. The authors thank
Nguyen~H.~V.~Hung as well as the other organizers of this
conference.

\section{Speculation concerning the localization of the double loop space of a sphere
at an odd prime $p$, and applications}\label{speculation}

The main goal of this section is to point out that if the
equivariant homology of either $M(t,k)$ or $M'(t,k)$ satisfies one
statement below, then these spaces provide a method for constructing
the localization at an odd prime $p$ of the double loop space of an
odd sphere. Some consequences are also given.

Let $R[\Sigma_k]$ denote the group ring of the symmetric group over
a commutative ring $R$ with $1$, and let $\mathcal S$ denote a left
$R[\Sigma_k]$-module. Let $X$ denote a path-connected Hausdorff
space with a free, right action of the symmetric group $\Sigma_k$.
Let $H_*(X/\Sigma_k; \mathcal S)$ denote the homology of the chain
complex $ C_*(X) \otimes_{\mathbb Z[\Sigma_k]} \mathcal S$ where $
C_*(X)$ denotes the singular chain complex of $X$.

Observe that the natural inclusion $$M(t,k) \to\ Conf(\mathbb C,k)
$$ induces a homomorphism $$H_*(M(t,k)/\Sigma_k; \mathcal S) \to\
H_*(Conf(\mathbb C,k)/\Sigma_k; \mathcal S).$$ If $t$ is equal to an
odd prime $p$, and $\mathcal S$ is the coefficient module given by
$\mathbb F_p(\pm 1)$ the field of $p$-elements with the action of
$\Sigma_k$ specified by the sign representation, then one question
is to decide whether this map induces an isomorphism in mod-$p$
homology. There is no strong evidence either way, although an
affirmative answer has interesting consequences which are described
below. The analogous question for $p = 2$ fails at once by Theorem
\ref{thm:no retraction}.

The reason for the interest in these particular homology groups is
the following observation implicit in \cite{C} as follows.

\begin{thm}
Let $\mathbb F$ denote a field. For each integer $i$ greater than
$0$, there is an isomorphism
$$ \oplus_{k \geq 0} H_{i-k(2n-1)}
(Conf(\mathbb C,k)/\Sigma_k, \mathbb F(\pm 1)) \to\
H_i(\Omega^2S^{2n+1}; \mathbb F).$$
\end{thm}

The next corollary follows at once.

\begin{cor}
Let $\mathbb F$ denote a field, and $p$ an odd prime.

\begin{enumerate}
    \item If the natural inclusion $$M(p,k) \to\ Conf(\mathbb C,k)$$
induces an isomorphism $$H_*(M(p,k)/\Sigma_k; \mathbb F_p(\pm 1))
\to H_*(Conf(\mathbb C,k)/\Sigma_k; \mathbb F_p(\pm 1)),$$ then
there are isomorphisms
$$ \oplus_{k \geq 0} H_{i-k(2n-1)}
(M(p,k)/\Sigma_k, \mathbb F_p(\pm 1)) \to\ H_i(\Omega^2S^{2n+1};
\mathbb F_p).$$
    \item If the natural inclusion $$M'(p,k) \to\ Conf(\mathbb R^2,k)$$
    induces an isomorphism
    $$H_*(M'(p,k)/\Sigma_k; \mathbb F_p(\pm 1)) \to\ H_*(Conf(\mathbb
    C,k)/\Sigma_k; \mathbb F_p(\pm 1)),$$ then there are isomorphisms
$$ \oplus_{k \geq 0} H_{i-k(2n-1)}
(M'(p,k)/\Sigma_k, \mathbb F_p(\pm 1)) \to\ H_i(\Omega^2S^{2n+1};
\mathbb F_p).$$
    \end{enumerate}
\end{cor}

The spaces $M(p,k)$ and $M'(p,k)$ are used next to give analogues of
labeled configuration spaces in which the configuration space itself
is replaced by a ``center of mass construction" as given above. Let
$Y$ denote a pointed space with base-point $*$ and $W$ any
topological space. Recall the labeled configuration space
$$C(W,Y)$$ given by equivalence classes of pairs $[S,f]$ where

\begin{enumerate}
    \item $S$ is a finite subset of $W$,
    \item $f: S \to Y$ is a function, and
    \item $[S,f]$ is equivalent to $[S-\{p\},f|_{S-\{p\}}]$ if and
    only if $f(p) = *$.
\end{enumerate}

One theorem proven in \cite{May} is as follows.
\begin{thm} \label{thm:May}
If $Y$ is a path-connected CW-complex, then $C(\mathbb R^n,Y)$ is
homotopy equivalent to $\Omega^n\Sigma^n(Y)$
\end{thm}

Technically, May's proof does not exhibit a map between these two
spaces. There are weak equivalences on the level of May's
construction \cite{May} $\alpha: C_n(Y) \to \Omega^n \Sigma^n(Y)$
and the natural evaluation map $e: C_n(Y) \to C(\mathbb R^n,Y)$.

Furthermore, the construction $D_k(\Omega^2\Sigma^2(X))$ is homotopy
equivalent to $$Conf(\mathbb C,k)
\times_{\Sigma_k}X^{(k)}/Conf(\mathbb C,k)\times_{\Sigma_k} \{*\}$$
for which $X^{(k)}$ denotes the $k$-fold smash product \cite{May}.
When localized at an odd prime $p$, $D_p(\Omega^2 S^{2n+1})$ is
homotopy equivalent to a mod-$p$ Moore space $P^{2np-1}(p)$ with a
single non-vanishing reduced homology group given by $\mathbb Z/
p\mathbb Z$ in dimension $2np-2$. This last assertion follows from
the computations in \cite{C}.

\begin{defn}
Define $$M_t(\mathbb C, Y)$$ to be the subspace of $C(\mathbb C, Y)$
given by those points for which $S$ is a subset of $M(t,k)$ with
natural inclusion denoted by $I_p:M_t(\mathbb C, Y)\to C(\mathbb C,
Y)$ and
$$M_t'(\mathbb C, Y)$$ to be the subspace of $C(\mathbb C, Y)$ given
by those points for which $S$ is a subset of $M'(t,k)$ with natural
inclusion denoted ( ambiguously ) by $I_p:M_t(\mathbb C, Y)\to
C(\mathbb C, Y)$.
\end{defn}

The next statement provides a potential method for constructing the
localization at $p$ of the double loop space of an odd sphere which
also has some useful properties.
\begin{thm} Assume that $p$ is an odd prime.

\begin{enumerate}
    \item If $M(t,k) \to\ Conf(\mathbb C,k) $ induces an isomorphism
$$H_*(M(t,k)/\Sigma_k; \mathbb F_p(\pm 1)) \to\ H_*(Conf(\mathbb C,k)/\Sigma_k; \mathbb F_p(\pm 1))$$ for $t$ an odd prime $p$, then
the natural map
$$I_p: M_p(\mathbb C, S^{2n-1}) \to\
\Omega^2S^{2n+1}$$ induces a mod-$p$ homology isomorphism.
    \item If $M'(t,k) \to\ Conf(\mathbb C,k) $ induces an isomorphism
$$H_*(M'(t,k)/\Sigma_k; \mathbb F_p(\pm 1)) \to\ H_*(Conf(\mathbb C,k)/\Sigma_k; \mathbb F_p(\pm 1))$$ for $t$ an odd prime $p$, then
the natural map $$I_p: M'_p(\mathbb C, S^{2n-1}) \to\
\Omega^2S^{2n+1}$$ induces a mod-$p$ homology isomorphism.
\end{enumerate}
\end{thm}

One consequence of this last theorem is that it implies properties
of the double suspension of $E^2: S^{2n-1}\to \Omega^2 S^{2n+1}$
after localization at an odd prime. In particular, the next
corollary follows directly.
\begin{cor}\label{cor:consequences}
Let $p$ denote an odd prime. If either (1) the natural inclusion
$M(p,k) \to\ Conf(\mathbb C,k)$ induces an isomorphism
$$H_*(M(p,k)/\Sigma_k; \mathbb F_p(\pm 1)) \to H_*(Conf(\mathbb C,k)/\Sigma_k; \mathbb F_p(\pm 1)),$$ or (2) the natural inclusion
$M'(p,k) \to\ Conf(\mathbb R^2,k)$ induces an isomorphism
$$H_*(M'(p,k)/\Sigma_k; \mathbb F_p(\pm 1)) \to\ H_*(Conf(\mathbb C,k)/\Sigma_k; \mathbb F_p(\pm 1)),$$ then after localization at
$p$, the mod-$p$ Moore space $$P^{2np+1}(p)$$ is a retract of
$\Sigma^2 \Omega^2 S^{2n+1}$. In that case, the following hold:
\begin{enumerate}
\item Any map  $$\alpha: P^{2p+1}(p) \to\ S^3$$
  given by an extension of $\alpha_1: S^{2p} \to S^3$, which
realizes the first element of order $p$ in the homotopy groups of
the $3$-sphere induces a split epimorphism on the $p$-primary
component of homotopy groups.
  \item After localization at the prime $p$, the homotopy theoretic
  fibre of the double suspension $E^2: S^{2n-1} \to \Omega^2 S^{2n+1}$ is the fibre of
 a map $\Omega^2 S^{2np+1}\to S^{2np-1}$.
\end{enumerate}
\end{cor}

{\bf Remarks:}
\begin{enumerate}
\item The main content of Corollary \ref{cor:consequences} is that the
``center of mass arrangements" may provide a useful way to construct
a localization at odd primes for the double loop space of an odd
sphere. Corollary \ref{cor:consequences} follows from Theorem
\ref{thm:comparison} as outlined in \cite{C2}.
\item In addition, Theorem \ref{thm:no retraction} shows that these
constructions fail to give the localization at the prime $2$ of
$\Omega^2 S^{2n+1}$.
\item It would also be interesting to see whether
there are analogous properties for the center of mass arrangement
with $\mathbb C$ replaced by $\mathbb C^n$ which may provide the
localization of $\Omega^{2n} S^{2(n+k)+1}$ at an odd prime $p$.
\end{enumerate}

\section{Sketch of Proposition \ref{prop:proposition cartesion}}
\label{proof of cartesion}

Notice that the set theoretic pull-back in the diagram given in
Proposition \ref{prop:proposition cartesion} is precisely the
subspace of the configuration space given by the $t$-fold center of
mass arrangement $M(t,k)$.

\section{Calculations at the prime 2, and the proof of Theorem \ref{thm:no retraction}}
\label{ the prime two }

The method here of comparing the homology of the center of mass
arrangement with that of the configuration space uses some
additional topology. Here, consider the natural inclusion $M(p,k)
\to\ Conf(\mathbb C,k)$ together with the induced map $$I_p:
M_p(\mathbb C, S^{2n-1}) \to\ C(\mathbb C, S^{2n-1}).$$

The space $C(\mathbb C, S^{2n-1})$ is homotopy equivalent to
$\Omega^2 S^{2n+1}$ \cite{May}. In addition, the spaces $M_p(\mathbb
C, S^{2n-1})$, and $C(\mathbb C, S^{2n-1})$ admit stable
decompositions which are compatible  by the remarks in
\cite{C2,CMT}. Notice that the inclusion $M(t,k)\to Conf(\mathbb
C,k)$ is the identity in case $k \leq t$ by definition of $M(t,k)$.
Thus the induced maps on stable summands
$$D_j(M_p(\mathbb C, S^{2n-1})) \to D_j(C(\mathbb C, S^{2n-1}))$$
is the identity in case $j \leq p,$ a feature which is used below.

Let $p = 2$, and consider the second stable summand
$$D_2(M_2(\mathbb C, S^{2n-1})) = D_2(C(\mathbb C, S^{2n-1})).$$
This stable summand is homotopy equivalent to
$$(S^1\times_{\Sigma_2}S^{4n-2})/ (S^1 \times _{\Sigma_2}*)$$ which is itself
homotopy equivalent to  $$\Sigma^{4n-3}(\mathbb R \mathbb P^2).$$
Let $u$ denote a basis element for $H_{4n-2}(\Sigma^{4n-3}(\mathbb
R\mathbb P^2);\mathbb F_2)$, and $v$ denote a basis element for
$H_{4n-1}(\Sigma^{4n-3}(\mathbb R\mathbb P^2);\mathbb F_2)$.

In addition, there is a strictly commutative diagram
\[
\begin{CD}
M_2(\mathbb C, S^{2n-1})  @>{h_2}>>  \Omega^{\infty}\Sigma^{\infty}(D_2(C(\mathbb C, S^{2n-1}))) \\
 @VV{I_2}V          @VV{1}V \\
C(\mathbb R^2, S^{2n-1}) @>{h_2}>>
\Omega^{\infty}\Sigma^{\infty}(D_2(C(\mathbb C, S^{2n-1})))
\end{CD}
\] which gives the fact that the space $D_2(C(\mathbb C, S^{2n-1})) = \Sigma^{4n-3}(\mathbb R\mathbb P^2)$
is a stable retract of both spaces, in a way which is compatible
with the natural stable decompositions.

Further, by Proposition \ref{prop:proposition cartesion}, together
with the definition \cite{C2} of the map $$h_2:M_2(\mathbb C,
S^{2n-1}) \to \Omega^{\infty}\Sigma^{\infty}(D_2(C(\mathbb C,
S^{2n-1}))),$$ there is a compression of this map through
$\Omega^{2}\Sigma^{2}(D_2(C(\mathbb C, S^{2n-1})))$. Thus, there is
a commutative diagram given as follows.

\[
\begin{CD}
M_2(\mathbb C, S^{2n-1})  @>{h_2}>>  \Omega^{2}\Sigma^{2}(D_2(C(\mathbb C, S^{2n-1}))) \\
 @VV{I_2}V          @VV{\Omega^2(E)}V \\
C(\mathbb R^2, S^{2n-1}) @>{h_2}>>
\Omega^{\infty}\Sigma^{\infty}(D_2(C(\mathbb C, S^{2n-1}))).
\end{CD}
\] These remarks have the following consequence for which $Q_i(x)$ is the standard
notation for Araki-Kudo-Dyer-Lashof operations as described in
\cite{C}.

\begin{lem} \label{lem: compressions}
The image of the map $ h_2: M_2(\mathbb C, S^{2n-1}) \to
\Omega^{\infty}\Sigma^{\infty}(D_2(M_2(\mathbb C, S^{2n-1}))) $ on
the level of mod-$2$ homology is contained in the subalgebra
generated by the elements $x$, and $Q_1^q(x)$ for $q \geq 1$ for
which $x$ is an element of a basis for the mod-$2$ homology of
$\Sigma^{4n-3}(\mathbb R\mathbb P^2)$ given by $\{u,v\}$. In
particular, the element $Q_3(x)$ cannot appear as a non-trivial
summand of the image.
\end{lem}

\begin{lem} \label{lem: surjection}
If $k \geq 4$, and the natural map
$$M(2,k)/\Sigma_k \to Conf(\mathbb C,k)/\Sigma_k $$
induces a surjection in mod-$2$ homology, then the natural map
$$M(2,4)/\Sigma_4 \to Conf(\mathbb C,4)/\Sigma_4 $$ induces a
surjection in mod-$2$ homology.
\end{lem}

\begin{proof}
Notice that if $k \geq 4$,  the space $Conf(\mathbb C,4)/\Sigma_4$
is a stable retract of the space $Conf(\mathbb C,k)/\Sigma_k$ via a
map induced by the transfer obtained from the natural
$\Sigma_k$-cover \cite{CMT2}. Thus there is a commutative diagram
\[
\begin{CD}
 \Sigma^{2k}(M(2,k)/\Sigma_k) @>{}>>  \Sigma^{2k}(Conf(\mathbb C,k)/\Sigma_k) \\
 @VV{tr}V          @VV{tr}V \\
\Sigma^{2k}(M(2,4)/\Sigma_4) @>{}>> \Sigma^{2k}(Conf(\mathbb
C,4)/\Sigma_4)
\end{CD}
\] in which the vertical maps are induced by the natural
transfer. Hence the natural map $M(2,4)/\Sigma_4 \to Conf(\mathbb
C,4)/\Sigma_4$ induces a surjection on mod-$2$ homology as the maps
$M(2,k)/\Sigma_k \to Conf(\mathbb C,k)/\Sigma_k$ as well as
$tr:\Sigma^{2k}(Conf(\mathbb C,k)/\Sigma_k) \to
\Sigma^{2k}(Conf(\mathbb C,4)/\Sigma_4)$ induce surjections on
mod-$2$ homology by \cite{CMT2}.

\end{proof}

The proof of Theorem \ref{thm:no retraction} is given next.
\begin{proof}

Assume that the natural inclusion $M(2,k) \to\ Conf(\mathbb C,k)$
induces an epimorphism on the level of $H_*(M(2,k)/\Sigma_k; \mathbb
F_2) \to\ H_*(Conf(\mathbb C,k)/\Sigma_k;\mathbb  F_2)$ for some $ k
\geq 4$. Then by Lemma \ref{lem: surjection}, the natural map
$M(2,4)/\Sigma_4 \to Conf(\mathbb C,4)/\Sigma_4$ induces a
surjection on mod-$2$ homology, and the induced map $H_*(M_2(\mathbb
C, S^{2n-1});\mathbb F_2) \to H_*(C(\mathbb C, S^{2n-1});\mathbb
F_2)$ is an epimorphism in dimensions $\leq 8n-1$. This will lead to
a contradiction.

Notice that
\begin{enumerate}
\item ${h_2}_*({x_{2n-1}}^2) = u$,

\item ${h_2}_*(Q_1(x_{2n-1})) = v$ and

\item ${h_2}_*(Q_1Q_1(x_{2n-1})) = AQ_1(v) + BQ_3(u)$ for scalars
$A$, and $B$ where $u$ is the unique non-zero class in
$H_{4n-2}(D_2(C(\mathbb R^2, S^{2n-1})); \mathbb F_2)$, and $v$ is
the unique non-zero class in $H_{4n-1}(D_2(C(\mathbb C, S^{2n-1}));
\mathbb F_2)$ \cite{C}.

\end{enumerate}

A direct computation using $Sq^1_*$, $Sq^2_*$, and the coproduct
gives $$A = B  =1.$$ The details are as follows. Notice that
$Sq^2_*(Q_1Q_1(x_{2n-1})) = 0$, but that $Sq^2_*(Q_1(v)) = Q_1(u) =
Sq^2_*(Q_3(u))$. Thus $A = B$. Furthermore $Sq^1_*(Q_1Q_1(x_{2n-1}))
= Q_1(x_{2n-1})^2.$

Finally notice that ${h_2}_*(({x_{2n-1}}^2)\cdot Q_1(x_{2n-1})) = u
\cdot v + P$ where $P$ is a primitive element. The only non-zero
choice for this primitive element $P$ is $Q_1(u)$. However,
$Sq^1_*(P) = 0$. Thus ${h_2}_*({x_{2n-1}}^4) = Sq^1_*(u \cdot v + P)
= u^2$. Hence
$$Sq^2_*Sq^1_*{h_2}_*(Q_1Q_1(x_{2n-1}))= u^2$$ and $A= B = 1$.

It follows that if the natural map $M(2,4)/\Sigma_4 \to Conf(\mathbb
C,4)/\Sigma_4$ induces a surjection on mod-$2$ homology, then the
class $Q_1(v) + Q_3(u)$ is in the image of the composite of the
following two maps:

$${I_2}_*: H_*(M_2(\mathbb C, S^{2n-1});\mathbb F_2) \to H_*(C(\mathbb C,
S^{2n-1});\mathbb F_2),$$ and

$${h_2}_*: H_*(C(\mathbb C, S^{2n-1});\mathbb F_2) \to
H_*(\Omega^{\infty}\Sigma^{\infty}D_2(C(\mathbb C,
S^{2n-1}));\mathbb F_2).$$

Thus the above computation gives that the class $Q_1(v) + Q_3(u)$ is
in the image of the composite
\[
\begin{CD}
H_*(M_2(\mathbb C, S^{2n-1});\mathbb F_2)  @>{{h_2}_* \circ
{I_2}_*}>> H_*(\Omega^{\infty}\Sigma^{\infty}D_2(C(\mathbb C,
S^{2n-1}));\mathbb F_2).
\end{CD}
\]

By Lemma \ref{lem: compressions}, the class $Q_1(v) + Q_3(u)$ cannot
be in the image of the map
$$H_*(\Omega^{2}\Sigma^{2}(D_2(M_2(\mathbb C, S^{2n-1})));\mathbb F_2)
\to H_*(\Omega^{\infty}\Sigma^{\infty}(D_2(C(\mathbb C,
S^{2n-1})));\mathbb F_2).$$ Hence, (1) the natural map
$M(2,4)/\Sigma_4 \to Conf(\mathbb C,4)/\Sigma_4$ cannot induce a
surjection on mod-$2$ homology and (2) the natural map
$M(2,k)/\Sigma_k \to Conf(\mathbb C,k)/\Sigma_k$ for $ k \geq 4$
cannot induce a surjection on mod-$2$ homology . The theorem
follows.
\end{proof}

\section{Sketch of Theorem \ref{thm: splitting} and Corollary \ref{cor:sign.rep.splitting}}
\label{proof of splitting}

The proof of follows Theorem \ref{thm: splitting} at once from the
constructions in the appendix of \cite{C2} or the main theorem in
\cite{CMT} where it was shown that these maps admit stable right
homotopy inverses.

To check Corollary \ref{cor:sign.rep.splitting}, notice that the
sign representation is given by the action of the symmetric group on
the top non-vanishing homology group of $(S^1)^n$. The corollary
follows from Theorem \ref{thm: splitting}.

\section{Sketch of Theorem \ref{thm:comparison}} \label{proof of comparison}
The commutativity of the diagram in \ref{thm:comparison} follows by
definition. That the map $$h_p: \Omega^2S^{2n+1} \to
\Omega^{2p}\Sigma^{2p}(D_p(\Omega^2\Sigma^2(S^{2n-1})))$$ induces an
isomorphism on $H_{2np-2}(-;\mathbb F_p)$ is checked in \cite{C2}.
Since $h_p$ induces an isomorphism on $H_{2np-2}(-;\mathbb F_p)$, it
follows from the known homology of these spaces that $\bar h_p$ does
also. Given a map with the homological properties of $\bar h_p$, the
proof of theorem \ref{thm:comparison} follows from \cite{C2}.

Remark: The goal of this approach is to try to desuspend a map
analogous to that given in \cite{selick}.

\bibliographystyle{amsalpha}

\end{document}